\documentclass[12pt, letterpaper, oneside, reqno]{amsart}
\usepackage{amsthm,amsfonts,amssymb,amsmath,eucal}
\newcommand{\la}{\lambda}

\newcommand{\Ga}{\Gamma}
\newcommand{\de}{\delta}
\newcommand{\ep}{\varepsilon}
\newcommand{\om}{\omega}

\newcommand{\Z}{\mathbb{Z}}

\newcommand{\sgn}{\operatorname{sgn}}
\newcommand{\pf}{\operatorname{Pf}}

\newcommand{\fd}{\downarrow}

\theoremstyle{plain}
\newtheorem{theor}{Theorem}[section]

\newtheorem{propo}[theor]{Proposition}
\newtheorem{corol}[theor]{Corollary}
\theoremstyle{definition}

\newtheorem{defin}[theor]{Definition}
\numberwithin{equation}{section}

\begin{document}

\title
{Interpolation analogues of
Schur Q-functions}
\author{Vladimir N.~Ivanov}
\address{Russia, 119992,
Moscow, Vorobjevy Gory,
Moscow State University, GZ,
Department of Mathematics and Me\-cha\-nics,
Chair of Higher Algebra}
\email{vivanov@mccme.ru}
\begin{abstract}
We introduce interpolation analogues of
Schur Q-func\-tions --- the multiparameter
Schur Q-functions. We obtain for them
several results:
a combinatorial formula, generating functions
for one-row and two-rows functions, vanishing
and characterization pro\-per\-ties, a Pieri-type
formula, a Nimmo-type
formula (a relation of two Pfaffians),
a Giambelli-Schur-type Pfaffian formula,
a determinantal formula for the transition
coefficients between multiparameter
Schur Q-functions with different parameters.
We write an explicit Pfaffian
expression for the dimension
of skew shifted Young diagram.
This paper is a continuation
of author's paper
math.CO/0303169 and is a partial projective analogue
of the paper by A. Okounkov and
G. Olshanski q-alg/9605042, and of the paper
by G. Olshanski, A. Regev and A. Vershik
math.CO/0110077.
\end{abstract}

\maketitle

\tableofcontents

\section{Introduction}

This text is an announcement of
the results. In most cases we give
only the main ideas of the proofs.
The detailed proofs will be presented
later.

In 1911 in the paper \cite{Sc} I.~Schur introduced
the Q-functions in order to study the projective
representations of the symmetric groups.
The projective representations of the symmetric group
$S(n)$ are linearized by the spin-symmetric group
$\widetilde S(n)$ (one of two possible
central $\Z_2$-extensions of the group $S(n)$).
The Schur Q-functions $Q_\la$, where $\la$ runs over
partitions of a number $n$ on distinct parts,
play the same role in the representation theory
of the spin-symmetric groups (\cite{Sc}, \cite{HH},
\cite{St1}, \cite{Jo})
as the conventional Schur S-functions in the
representation theory of the symmetric groups
(\cite[Ch.I,\S7]{Ma1}). A "projective"
analogue of the Schur-Weyl duality for $gl(N)$ and
$S(n)$ (\cite{We}) was obtained by A.~Sergeev in
\cite{Se1}. A.~Sergeev introduced a new group
(now it is called the Sergeev group) instead of $S(n)$
and replaced the Lie algebra $gl(N)$ by the Lie
superalgebra $Q(N)$. Thus the Schur Q-functions,
up to scalar factors, are equal to the
characters of the irreducible
representations of the Lie superalgebra $Q(N)$,
see also \cite{Ya}.
Besides the Schur Q-functions it is convinient to use
Schur P-functions $P_\la$, which differ from
the Q-functions $Q_\la$ by simple scalar
factors (Definition~\ref{d-qfun}).
The Schur P-functions
are $t=-1$ specializations
of  the Hall-Littlewood polynomials
and constitute a linear basis
in the algebra of supersymmetric functions $\Gamma$,
see \cite[\S8]{Ma1}, \cite{Pr}, \cite{HH} for
more details.
For a partition with distinct parts
$\la=(\la_1>\la_2>\dots>\la_l>0)$ the Schur P-function
$P_\la$ is the symmetric polynomial (suppose
the number of variables
$n$ is greater than $l$)
\begin{gather}
P_{\la}(x_1,\dots,x_n)
=\frac1{(n-l)!}
\sum_{\om\in S(n)}
\prod_{i=1}^l x_{\om(i)}^{\la_i}
\prod_{i\le l,\,i\le j\le n}\frac
{x_{\om(i)}+x_{\om(j)}}
{x_{\om(i)}-x_{\om(j)}};\\
Q_{\la}(x_1,\dots,x_n)=
2^l P_{\la}(x_1,\dots,x_n),
\end{gather}
see Definition~\ref{d-trad}
for the details.

We replace the ordinary powers $x^k$ by
the generalized powers $(x-a_1)(x-a_2)\dots (x-a_k)$
(we suppose $a_1=0$) in these definitions of
the Schur P- and Q-functions and we obtain
(Definition~\ref{d-mult}) new supersymmetric functions
$P_{\la;a}$ and
$Q_{\la;a}$.
We call them the multiparameter Schur P- and Q-functions.
The classical Schur Q-function $Q_\la$ is the
leading homogenous term of the function $Q_{\la;a}$.
As the ordinary Schur Q-functions
$Q_\la$ the multiparameter Schur Q-functions
$Q_{\la;a}$ constitute a linear basis in
the algebra of supersymmetric functions.
When all $a_j$ are distinct, the function
$Q_{\la;a}$ can be viewed as a solution
of some multivariate interpolation problem,
see Section~\ref{s-char-prop}. For this reason
we call
these functions the interpolation analogues of
the Schur Q-functions.

Earlier
interpolation analogues for many
other symmetric functions were found and
studied: for analogues of the classical Schur
S-functions, see \cite{BL},
\cite{Ok1}, \cite{OO1}, \cite{OO2}, and references
in \cite{OO1};
in the case of the supersymmetric
Schur S-functions, see \cite{Mo},
\cite{ORV}; in the case of the Jack functions,
see
\cite{KS}, \cite{KOO}, \cite{OO3};
the case of the Macdonald polynomials, excepting
the Hall-Littlewood polynomials, was considered in
\cite{Kn}, \cite{Sa}, \cite{Ok2}, \cite{Ok3},
and references therein.

We obtain the following main results about
the multiparameter Schur Q-functions:
\begin{itemize}
\item A Nimmo-type formula
expressing $P_{\la;a}$ as a relation
of two Pfaffians (Section~ \ref{s-Nimm}).
\item A combinatorial formula expressing
$Q_{\la;a}$ in terms of marked shifted tableaux
(Section~\ref{s-comb}).
This formula has an additional symmetry
as compared to its counterpart for the
conventional
Schur Q-functions. Also we rewrite our formula
in terms of unmarked shifted tableaux.
\item Generating functions for the
one-row $Q_{(r);a)}$ and the two-rows
$P_{(r,s);a}$ multiparameter functions
(Section~\ref{s-gene}). (Note that
even in the ordinary case the formula
for the generating function for the
two-rows $P_{(r,s)}$ seems to be new)
\item A Giambelli-Schur-type Pfaffian
formula that expresses an arbitrary
multiparameter Schur Q-function
$Q_{\la;a}$ as a Pfaffian
of the two-rows multiparameter
Schur Q-functions (Section~\ref{s-Giam}).
\end{itemize}
Also we find explicitly the transition
coefficients between two bases
$\{Q_{\la;a}\}$ and $\{Q_{\la;b}\}$
for two sequences $(a_k)$ and
$(b_k)$ (Section~\ref{s-tran}).
In particular we may
express each of the multiparameter
Schur Q-functions $Q_{\la;a}$
as a linear combination of the ordinary
Schur Q-functions $Q_\mu$.

Suppose $a_k=k-1$; then we obtain
an important particular case of
the multiparameter Schur Q-functions.
We call them the factorial Schur Q-functions.
In Section~\ref{s-dime} we calculate
the number of shifted standard
tableaux of a given skew shape
in terms of the factorial
P-functions
(called a
dimension of a skew shifted
Young diagram).
We obtain for this dimension
a simple Pfaffian expression
(Theorem~\ref{t-dime-pfaf}).
This dimension
may be also
rewritten in terms of
the multiplicity of the restriction of
irreducible representations
of the spin-symmetric group
or the Sergeev group to the smaller
subgroup, for the details see,
for example, \cite{HH}, \cite{Iv3}.
Using this formula with the factorial
P-functions in \cite{Iv1}
we obtain a new proof of
Nazarov's theorems about
the characters of the infinite
spin-symmetric group (\cite{Na1}).
Note also that the factorial Schur P-functions
possibly may be obtained from the Schubert
polynomials (\cite{La1}, \cite{La2}, \cite{LP},
\cite{FK}, \cite{BH}) but
we did not verify this fact.
Certain problems related to
the factorial Schur Q-functions are
considered in \cite{Na2}, \cite{Iv3}.

The author is very grateful to G.~Olshanski
for his constant attention to this work
and many valuable remarks, to A.~Okounkov
for the main Definition~\ref{d-fact},
to A.~Borodin for the statement of
Theorem~\ref{t-dvus} (The generating function
for the two-rows multiparameter P-functions).

\section{Notation and definitions}
\label{s-nota}

\begin{defin}
A polynomial $f(x_1,\dots,x_n)$ is
called {\it supersymmetric\/} if
the following conditions hold:
\begin{enumerate}
\item $f(x_1,\dots,x_n)$ is symmetric in
variables $x_1,\dots,x_n$ ;
\item for all $1\le i<j\le n$ the polynomial
$$
f(x_1,\dots,x_{i-1},t,x_{i+1},\dots,x_{j-1},-t,x_{j+1},\dots,x_n)
$$
does not depend on $t$ .
\end{enumerate}
\end{defin}

Supersymmetric polynomials in $n$ variables
form an algebra .
Denote this algebra
by $\Ga(n)$ .
The algebra $\Ga(n)$ is graded
by degree of polynomials. The specialization
$x_{n+1}=0$ is a morphism of
the graded algebras
\begin{equation} \label{proj}
\Ga(n+1)\to\Ga(n)\,.
\end{equation}

\begin{defin} \label{defGa}
\emph{The algebra $\Ga$ of supersymmetric functions}
is the projective limit
$$
\Ga=\varprojlim\Ga(n), \quad n\to\infty,
$$
in the category of graded algebras, taken with respect to morphisms
(\ref{proj}).
In other words,
an element $f\in\Ga$ is a
sequence $(f_n)_{n\ge1}$ such that:
\begin{enumerate}
\item $f_n\in\Ga(n)$, $n=1,2,\dots$,
\item $f_{n+1}(x_1,\dots,x_n,0)=f_{n}(x_1,\dots,x_n)$
(the stability condition),
\item $\sup_n\deg f_n<\infty$\,.
\end{enumerate}
\end{defin}

\begin{defin}\label{d-part}
A finite sequence $\la=(\la_1\ge\ldots\ge\la_k\ge 0)$
of decreasing positive integers is called
\emph{a partition}.
The number of non-zero parts of $\la$
is called \emph{the length of the partition} $\la$
and is denoted by $\ell(\la)$.
Let us denote by $|\la|$ the sum of the parts
of the partition $\la$:
$$
|\la|=\sum_{k=1}^{\ell(\la)}
\la_k.
$$
If $|\la|=n$ we shall also
write $\la\vdash n$.
Also put $m_i(\la)=\#\{k\mid\la_k=i\}$ for
$i\ge 1$. Obviously, we have
$$
\ell(\la)=\sum_{i\ge 1} m_i(\la),\quad
|\la|=\sum_{i\ge 1} i m_i(\la).
$$
\end{defin}

\begin{defin}\label{d-stri}
A partition is called {\it strict\/} if  any two
non-zero parts of it are distinct. We shall
denote by $DP_n$ the set of the strict partitions
of the number $n$.
A partition is called {\it odd\/} if  all
non-zero parts of it are odd. We shall
denote by $OP_n$ the set of the odd partitions
of the number $n$. Also put
\begin{equation}
DP=\bigcup_{k\ge 0} DP_k,\quad
OP=\bigcup_{k\ge 0} OP_k.
\end{equation}
\end{defin}

\begin{defin}\label{d-shif-Youn}
For the strict partition one may define
the shifted Young diagram besides the simple
Young diagram
\cite[Ch.I, \S1]{Ma1}. If
$\la$
is a strict partition, then the set
$$
\{(i,j)\in \Z^2\,|\, i\le j\le \la_i+i-1, 1\le i\le \ell(\la)\}
$$
is called {\it the shifted Young diagram of the partition
$\la$}
and is denoted by
$D(\la)'$.
It is useful to draw unit squares
instead of points of $\Z^2$.
We assume that the first coordinate axis is directed
downwards, the second coordinate axis is directed
to the right and the point $(0,0)$ is in the left
upper corner of the figure.
For example,
if $\la=(7,4,3,1)$, then $D(\la)'=$
\bigskip
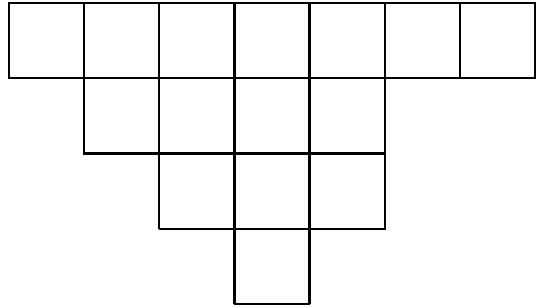
\begin{center}
\unitlength=0.5cm
\begin{picture}(14,8)
\put(0,8){\line(1,0){14}}
\put(0,6){\line(1,0){14}}
\put(2,4){\line(1,0){8}}
\put(4,2){\line(1,0){6}}
\put(6,0){\line(1,0){2}}
\put(0,8){\line(0,-1){2}}
\put(2,8){\line(0,-1){4}}
\put(4,8){\line(0,-1){6}}
\put(6,8){\line(0,-1){8}}
\put(8,8){\line(0,-1){8}}
\put(10,8){\line(0,-1){6}}
\put(12,8){\line(0,-1){2}}
\put(14,8){\line(0,-1){2}}
\end{picture}
\end{center}
\end{defin}
\bigskip

\begin{propo}\label{p-vspo} Let $r(x_1,\dots,x_l)$ be
a polynomial in variables $x_1,\dots,x_l$.
For arbitrary
$n\ge l$ we put
\begin{equation}\label{def-of-r}
R_n(x_1,\dots,x_n)=r(x_1,\dots,x_l)
\prod\limits_{i\le l,i<j\le n}
\frac{x_i+x_j}{x_i-x_j}
\end{equation}
and
\begin{equation}\label{sum-r}
\widetilde R_n(x_1,\dots,x_n)=\sum\limits_{\om\in S(n)}
R_n(x_{\om(1)},\dots,x_{\om(n)}).
\end{equation}
Then we have
\begin{itemize}
\item[a)] $\widetilde R_n$ is the polynomial in $x_1,\dots,x_n$ and
$$
\deg\widetilde R_n\le\deg r;
$$
\item[b)] $\widetilde R_n$ is the supersymmetric polynomial;
\item[c)] \label{e-symm}
if for some $i\ne j$ the polynomial
$r(x_1,\dots,x_l)$ is symmetric in $x_i,x_j$,
then
$\widetilde R_n=0$;
\item[d)] if $r(x_1,\dots,x_l)$ is divisible by
$x_1 x_2\dots x_l$,
then
$$
\widetilde R_{n+1}(x_1,\dots,x_n,0)=(n+1-l)\widetilde R_n(x_1,\dots,
x_n).
$$
\end{itemize}
\end{propo}

\begin{proof} See \cite[Proposition~1.1]{Iv1}.
\end{proof}

Next we consider three special cases
of Proposition~\ref{p-vspo}
(Definitions~\ref{d-trad}, \ref{d-fact}, and
\ref{d-mult}).
First let us introduce generalizations
of the ordinary power $x^n$.

\begin{defin}\label{d-degr}
Put
$$
(x\fd n)=\prod_{k=1}^n (x-k+1)
$$
for $n\ge 1$.
Also put $(x\fd 0)=1$.
Suppose $(a_k)_{k\ge 1}$
is an arbitrary sequence of
complex numbers. Put
$$
(x\mid a)^n=\prod_{k=1}^n (x-a_k),\quad
(x\mid a)^0=1.
$$

\end{defin}

\begin{defin}\label{d-trad}
Suppose $\la$ is a partition, $\ell(\la)=l<n$.
Put
$$
r(x_1,\dots,x_l)=\frac{\prod^l_{i=1}x^{\la_i}_i.}
{(n-l)!}
$$
in Proposition~\ref{p-vspo}.
Namely, put
$$
P_{\la|n}(x_1,\dots,x_n)=\frac1{(n-l)!}
\sum_{\om\in S(n)}
\prod_{i=1}^l x_{\om(i)}^{\la_i}
\prod_{i\le l,\,i\le j\le n}\frac
{x_{\om(i)}+x_{\om(j)}}
{x_{\om(i)}-x_{\om(j)}}.
$$
If $\la\in DP$ then the polynomials $P_{\la|n}$
are the specializations of the Hall-Littlewood
polynomials when the parameter $t=-1$,
see \cite[Ch.III, \S2, \S8]{Ma1}.
If $\ell(\la)>n$ then we put $P_{\la|n}=0$.
From Proposition~\ref{p-vspo}
it follows that the sequence
$(P_{\la|n})_{n\ge1}$  defines
the supersymmetric function,
which is  denoted by $P_{\la}$.
The functions $P_{\la}$ are called
{\it the Schur P-functions\/}.
Note that if $\nu$ is not a strict partition then
from Proposition~\ref{p-vspo}c)
it follows that
$P_\nu=0$. So our definition differs at this point
from the traditional
definition of the Schur P-functions,
see \cite[Ch.III, \S2, (2.2)]{Ma1}.
\end{defin}

Next definition is due to A.~Okounkov.

\begin{defin}\label{d-fact}
Suppose
$l=\ell(\la)\le n,\,
\la\in DP$.
Put
$$
F_{\la|n}(x_1,\dots,x_n)=\prod\limits^l_{i=1}
(x_i\fd\la_i)
\prod\limits_{i\le l,i<j\le n}
\frac{x_i+x_j}{x_i-x_j}.
$$
By definition, we put
$$
P^*_{\la|n}=\frac1{(n-l)!}\sum\limits_{\om\in S(n)}
F_{\la|n}(x_{\om(1)},\dots,x_{\om(n)}).
$$
It corresponds to the case
$$
r(x_1,\dots, x_l)=\frac{\prod_{i=1}^l (x_i\fd \la_i)}
{(n-l)!}
$$
in Proposition~\ref{p-vspo}.
If $\ell(\la)>n$ then, by definition, put $P^*_{\la|n}=0$.
From Proposition~\ref{p-vspo} it follows that
the sequence
$(P^*_{\la|n})_{n\ge1}$ defines the supersymmetric function
$P^*_\la$.
We call these functions
{\it the factorial Schur P-functions}.
\end{defin}

Next definition generalizes both Definition~\ref{d-trad}
and Definition~\ref{d-fact}.

\begin{defin}\label{d-mult}
Suppose $a=(a_k)_{k\ge 1}$ is an arbitrary sequence
of complex numbers, $\la\in DP,\,\ell(\la)=l< n$.
Then put
$$
P_{\la;a|n}(x_1,\dots,x_n)=\frac1{(n-l)!}
\sum_{\om\in S(n)}
\prod_{i=1}^l (x_{\om(i)}\mid a)^{\la_i}
\prod_{i\le l,\,i\le j\le n}\frac
{x_{\om(i)}+x_{\om(j)}}
{x_{\om(i)}-x_{\om(j)}}.
$$
If $\ell(\la)>n$ then we put $P_{\la;a|n}=0$.
From Proposition~\ref{p-vspo}
it follows that if $a_1=0$ then the sequence
$(P_{\la;a|n})_{n\ge1}$  defines
the supersymmetric function,
which is  denoted by $P_{\la;a}$.
The functions $P_{\la;a}$ are called
{\it the multiparameter Schur P-functions\/}.
Further in the text we suppose $a_1=0$.
\end{defin}

\begin{propo} \label{p-star} Suppose $\la\in DP;\,
\ell(\la)\le n$, then
$$
P_{\la;a}(x_1,\dots,x_n)=P_{\la}(x_1,\dots,x_n)+
g(x_1,\dots,x_n),
$$
where $g(x_1,\dots,x_n)$ is a supersymmetric
polynomial such that $\deg g <|\la|$.
\end{propo}

\begin{proof}
It follows from Proposition~\ref{p-vspo},
Definition~\ref{d-trad} and Definition~\ref{d-mult}.
\end{proof}

\begin{propo}\label{p-basi}
\begin{enumerate}
\item[a)] \label{e-clas}
The set $\{P_{\la}\mid\la\in DP\}$ is
a linear basis of the algebra $\Gamma$.
\item[b)] The set $\{P_{\la;a}\mid\la\in DP\}$ is
a linear basis of the algebra $\Gamma$.
In particular, the set
$\{P_{\la}^*\mid\la\in DP\}$ is
a linear basis of the algebra $\Gamma$.
\end{enumerate}
\end{propo}

\begin{proof}
\begin{enumerate}
\item[a)] It follows from
\cite[Theorem 2.11]{Pr}.
\item[b)] It follows from the assertion
a) and Proposition~\ref{p-star}.
\end{enumerate}
\end{proof}

Suppose the partition $\nu\notin DP$;
then from Proposition~\ref{p-vspo} it follows that
$P_{\nu;a}\equiv P_\nu^*\equiv 0$.

\begin{defin}\label{d-qfun}
For arbitrary partition $\la$ we put
$$
Q_\la=2^{\ell(\la)} P_\la,\quad
Q^*_\la=2^{\ell(\la)} P^*_\la,\quad
Q_{\la; a}=2^{\ell(\la)} P_{\la; a}.
$$
The supersymmetric functions $Q_\la$
were introduced by I.~Schur in
\cite{Sc}. They are called
\emph{the Schur Q-functions}.
The supersymmetric functions $Q^*_\la$
($Q_{\la; a}$) are called
\emph{the factorial (multiparameter)
Schur Q-functions}.
\end{defin}

\section{Nimmo-type formula}\label{s-Nimm}

In this section we
obtain the formula for $P_{\la;a}$,
which is an analogue
of the formula for the ordinary Schur P-functions
obtained by Nimmo (\cite{Ni}).

Recall the definition of the Pfaffian
of a skew-symmetric matrix.

\begin{defin}\label{d-pfaf}
Suppose $A$ is a skew-symmetric matrix $2n\times 2n$. Put
$$
\pf(A)=\sum_{\om\in \tilde S(2n)}
\sgn(\om)\prod^n_{i=1}a_{\om(2i-1)\om(2i)},
$$
where the sum is over $\om\in\tilde S(2n)\subset
S(2n)$ such that
$$
\om(2i-1)<\om(2i)\quad\text{and}\quad
\om(1)<\om(3)<\dots<\om(2n-3)<\om(2n-1).
$$
\end{defin}

\begin{theor}\label{t-Nimm}
Suppose $n\ge \ell(\la)=l,\:\la\in DP,\:
a=(a_k)_{k\ge 1}$ is an arbitrary
sequence of complex numbers, $a_1=0$.
Let $A_0(x_1,\dots,x_n)$
denote the skew-symmetric $n\times n$
matrix
$$
\left(\frac{x_i-x_j}
{x_i+x_j}\right)_{1\le i, j\le n}
$$
and let $B_\la$ denote the $n\times l$ matrix
$$
\left((x_i | a)^{\la_{l+1-j}}
\right)_{i\le n,\, j\le l}.
$$
Let $A_\la (x_1,\dots, x_n)$ be the skew-symmetric
$(n+l)\times (n+l)$ matrix
$$
\left(
\begin{array}{cc}
A_0(x_1,\dots,x_n)&
B_\la(x_1,\dots,x_n)\\
-B_\la(x_1,\dots,x_n)^t&
0
\end{array}
\right).
$$
Put
$$
\pf_0(x_1,\dots,x_n)=
\left\{
\begin{array}{ll}
\pf(A_0(x_1,\dots, x_n))&
\mbox{ if $n$ is even};\\
\pf(A_0(x_1,\dots, x_n,0))&
\mbox{ if $n$ is odd}.
\end{array}
\right.
$$
Also put
$$
\pf_\la(x_1,\dots,x_n)=
\left\{
\begin{array}{ll}
\pf(A_\la(x_1,\dots, x_n))&
\mbox{ if $n+l$ is even};\\
\pf(A_\la(x_1,\dots, x_n,0))&
\mbox{ if $n+l$ is odd}.
\end{array}
\right.
$$
Then
$$
P_{\la;a}(x_1,\dots, x_n)=
\frac
{\pf_\la(x_1,\dots,x_n)}
{\pf_0(x_1,\dots,x_n)}.
$$
\end{theor}

\begin{proof}
The proof follows Nimmo's method,
see \cite{Ni} or
\cite[Ch.III, \S8, Example 13]{Ma1}.
\end{proof}

\section{Combinatorial formulas}\label{s-comb}

\begin{defin}\label{d-alfa}
Consider the ordered alphabet
$\mathbb{P}_n=\{1'<1<2'<2<\dots<n'<n\}$.
By definition, put
$$
|k'|=|k|=k,\, \sgn(k)=-\sgn(k')=1
$$
for an arbitrary natural number $k$.
\end{defin}

\begin{defin}
Let $\la$ be an arbitrary strict partition.
{\it A marked shifted Young tableaux of shape
$\la$ and length n} is a mapping
$T:D'_{\la}\to \mathbb{P}_n$
such that the following conditions hold:
\begin{enumerate}
\item $T(i,j)\le T(i+1,j),\,T(i,j)\le T(i,j+1)$;
\item for each natural number $k$
there is at most one $k$ in the image
of each column and at most one $k'$
in the image of each row of $D'_{\la}$.
\end{enumerate}
Let us denote by $MSTab(\la,n)$ the set
of all shifted marked Young tableaux of
shape $\la$ and length $n$.
\end{defin}

Here is an example of the shifted Young tableaux
of shape $(7,4,3,1)$ and length $6$:
\bigskip
\begin{center}
\unitlength=0.5cm
\begin{picture}(14,8)
\put(0,8){\line(1,0){14}}
\put(0,6){\line(1,0){14}}
\put(2,4){\line(1,0){8}}
\put(4,2){\line(1,0){6}}
\put(6,0){\line(1,0){2}}
\put(0,8){\line(0,-1){2}}
\put(2,8){\line(0,-1){4}}
\put(4,8){\line(0,-1){6}}
\put(6,8){\line(0,-1){8}}
\put(8,8){\line(0,-1){8}}
\put(10,8){\line(0,-1){6}}
\put(12,8){\line(0,-1){2}}
\put(14,8){\line(0,-1){2}}
\put(0.7,6.7){$1$}
\put(2.7,6.7){$1$}
\put(4.7,6.7){$2$}
\put(6.7,6.7){$3'$}
\put(8.7,6.7){$4'$}
\put(10.7,6.7){$4$}
\put(12.7,6.7){$6'$}
\put(2.7,4.7){$2'$}
\put(4.7,4.7){$3'$}
\put(6.7,4.7){$3$}
\put(8.7,4.7){$4'$}
\put(4.7,2.7){$3$}
\put(6.7,2.7){$4'$}
\put(8.7,2.7){$4$}
\put(6.7,0.7){$4'$}
\end{picture}
\end{center}
\bigskip

\begin{theor}(The combinatorial
formula for the multiparameter
Schur Q-functions)\label{t-comb}
$$
Q_{\la; a}(x_1,x_2,\dots,x_n)=
\sum_{T\in MSTab(\la,n)}\prod_{(i,j)\in D'_{\la}}
(x_{|T(i,j)|}-\sgn(T(i,j))a_{j-i+1}).
$$
\end{theor}

\begin{proof}
The combinatorial formula for the
factorial Schur Q-functions
(see Corollary~\ref{c-comb} below) is proved in \cite{Iv2}.
We may use the same method in the case of
the multiparameter Schur Q-functions with evident changes.
\end{proof}

\begin{corol} (The combinatorial
formula for the factorial
Schur Q-fun\-c\-tions)\label{c-comb}
$$
Q^*_{\la}(x_1,x_2,\dots,x_n)=
\sum_{T\in MSTab(\la,n)}\prod_{(i,j)\in D'_{\la}}
(x_{|T(i,j)|}-\sgn(T(i,j))(j-i)).
$$
\end{corol}

The combinatorial formulas for the other interpolation
analogues of the symmetric functions were obtained
earlier in \cite[\S11]{OO1}, \cite{Ok2}, \cite{Mo},
\cite{ORV}.

We may reformulate Theorem~\ref{t-comb} in terms of
unmarked shifted tab\-leaux.

\begin{defin}\label{d-bord}
Suppose $\mu$ and $\la$ are strict partitions,
$D(\mu)'\subseteq D(\la)'$; then we write
$\mu\subset\la$.
Let us consider in this case
\emph{the skew shifted diagram}
$D(\la)'\setminus D(\mu)'$.
For each box $b\in D(\la)'\setminus D(\mu)'$
let us denote by $i(b)$ and $j(b)$ its first
and second coordinates respectively.
The skew shifted diagram
$\nu=D(\la)'\setminus D(\mu)'$ is called
\emph{a border strip} if the following conditions
hold:
\begin{enumerate}
\item
The diagram
$\nu$ has no $2\times 2$ block of squares.
\item \label{numer-ii}
The set $I(\nu)=\{j(b)-i(b)\mid b\in\nu\}
\subset\Z_+\cup\{0\}$ is an interval
of integers, i.~e.
$I(\nu)=[\min I(\nu),\max I(\nu)]\cap
\Z$. When this condition is
satisfied we say that $\nu$ is
\emph{connected}.
\end{enumerate}

Following \cite[\S4]{ORV}
consider \emph{interior sides}
of the squares of the shape $\nu$:
\emph{an interior side} is adjacent
to two squares of $\nu$;
the total number of the interior sides
is equal to $|\nu|-1$.
To each interior side $s$ we attach
the coordinates $(\ep ,\de)$ of its
midpoint and we write $s=(\ep,\de)$.
Note that one of the coordinates
is always half-integer while
another coordinate is integer.

Suppose $a=(a_n)_{n\ge 1}$ is a sequence
of complex numbers, $\nu$ is a border strip;
then put
$$
f(\nu; a; x)= 2x
\prod'(x+a_{\de-\ep+1/2})
\prod'' (x-a_{\de-\ep+1/2}) ,
$$
where $\prod'$
is taken over
the horisontal interior sides $s=(\ep,\de)$
of $\nu$,
and $\prod''$
is taken over
the vertical interior sides $s=(\ep,\de)$
of $\nu$.

Suppose $\nu=D(\la)'\setminus D(\mu)'$
has no $2\times 2$ block of squares;
then we may represent $\nu$ as a disjoint
union of the minimal number of border stripes
$\nu=\cup_{j=1}^k \nu_j$, i.~e. such that
$\forall i,j$ the set $\nu_i\cup\nu_j$ is not
connected in the sence of the condition
\ref{numer-ii}. In this case we put
$$
f(\nu; a; x)=
\prod_{j=1}^k
f(\nu_j; a; x).
$$

\emph{A (unmarked) shifted tableau} $T$ of
shape $\la$ and length $n$ is a
sequence of strict partitions
$\varnothing=\la^{(0)}
\subset\la^{(1)}\subset\dots
\la^{(n)}=\la$ such that each
$D(\la^{(j)})'\setminus
D(\la^{(j-1)})'$ has no $2\times 2$
block of squares. In this case we put
$$
f(T; a; x_1,\dots, x_n)=
\prod_{j=1}^n
f((D(\la^{(j)})'\setminus
D(\la^{(j-1)})');
a; x_j).
$$
\end{defin}

\begin{theor}
As always we suppose that $a=(a_n)_{n\ge 1}$
is a sequence of complex numbers such that
$a_1=0$; then we have
$$
Q_{\la;a}(x_1,\dots, x_n)=
\sum_T f(T;a;x_1,\dots,x_n)
$$
summed over all shifted tableaux of shape
$\la$ and length $n$.
\end{theor}

\begin{proof}
We may easily deduce this theorem
from Theorem~\ref{t-comb}.
\end{proof}

\section{Characterisation properties}\label{s-char-prop}

\begin{defin}
Suppose $\mu$ is a strict partition,
$a=(a_n)_{n\ge 1}$ is a sequence of complex
numbers; then,
by definition, we put
\begin{multline}
H_a(\mu)=
\prod_{k=1}^{\ell(\mu)}
(a_{\mu_k+1}\mid a)^{\mu_k}
\prod_{i<j} \frac
{a_{\mu_i+1}+a_{\mu_j+1}}
{a_{\mu_i+1}-a_{\mu_j+1}}=\\
\prod_{i<j}\left((a_{\mu_i+1}+a_{\mu_j+1})
\prod_{k=\mu_j+2}^{\mu_i}(a_{\mu_i+1}-a_k)\right).
\end{multline}
\end{defin}

In the "factorial" case, i.~e., when $a_k=k-1$
the expression $H_a(\mu)$ becomes the "shifted
product of hook lengths"
(see \cite[Ch.III, \S7, ex\-am\-p\-les]{Ma1}):
$$
H(\mu)=\prod^{\ell(\mu)}_{t=1}
\mu_t!\prod_{i<j}\frac{\mu_i+\mu_j}
{\mu_i-\mu_j}.
$$

\begin{defin}
Suppose $\la$ is a strict partition,
$a=(a_n)_{n\ge 1}$ is a sequence of
complex numbers; then let us define the collection
of variables $x(\la)$ by means of formula
$$
x(\la)_i=a_{\la_i+1}.
$$
\end{defin}

\begin{theor}(Vanishing property)\label{t-vani}
Assume that the numbers $a_j$ are
pairwise distinct and $a_1=0$.
\begin{itemize}
\item[a)] If $\mu,\la\in DP,\:
\mu\not\subset\la$ then
$Q_{\mu;a}(x(\la))=P_{\mu;a}(x(\la))=0.$
\item[b)] $P_{\mu;a}(x(\mu))=H_a(\mu).$
\end{itemize}
\end{theor}

\begin{proof}
\begin{itemize}
\item[a)] This claim follows from
Definition~\ref{d-mult}.
\item[b)] This formula can be readily obtained by
a direct computation.
\end{itemize}
\end{proof}

Next we write explicitly Theorem~\ref{t-vani}
in the "factorial" case.

\begin{corol}
\begin{itemize}
\item[a)] If $\mu,\la\in DP,\:
\mu\not\subset\la$ then
$Q^*_{\mu}(\la_1,\dots,\la_{\ell(\la)})=
P_\mu^*(\la_1,\dots,\la_{\ell(\la)})=0.$
\item[b)] $P^*_\mu(\mu)=H(\mu).$
\end{itemize}
\end{corol}

Following the method of
\cite{Ok1}, \cite{OO1} we
may obtain characterization theorems
for the multiparameter Schur P-functions.

\begin{theor} (Characterization Theorem I)
\label{t-char-I}
Assume that the numbers $a_j$ are pairwise
distinct. Suppose $f\in \Ga$ satisfies the
following conditions:
\begin{enumerate}
\item The top degree homogenous component of
$f$ is equal to
$$
\sum_{\mu\in DP_n} c_\mu P_\mu.
$$
\item $\forall\la\in DP$ such that
$|\la|<n$ we have $f(x(\la))=0.$
\end{enumerate}
Then
$$
f\equiv \sum_{\mu\in DP_n}
c_\mu P_{\mu;a}.
$$
\end{theor}

\begin{proof}
We may use Proposition~\ref{p-star}
and Theorem~\ref{t-vani}.
\end{proof}

\begin{theor} (Characterization Theorem II)
\label{t-char-II}
Assume that the numbers $a_j$ are pairwise distinct.
Suppose $f\in\Ga$ satisfies the following conditions
for some strict partition $\mu$:
\begin{enumerate}
\item $f(x(\mu))=H_a(\mu).$
\item $\deg f\le |\mu|.$
\item $\forall\la\in DP$ such that
$|\la|\le |\mu|$ and $\la\ne\mu$
we have $f(x(\la))=0$.
\end{enumerate}
Then $f\equiv P_{\mu; a}$.
\end{theor}

\begin{proof}
Again we use Proposition~\ref{p-star}
and Theorem~\ref{t-vani}.
\end{proof}

Note that these characterization theorems
may be easy reformulated in the particular
"factorial" case.

Theorem~\ref{t-char-I} and
Theorem~\ref{t-char-II}
show that
the multiparameter Schur P- and Q-functions
can be considered as the
solutions of the multivariate interpolation
problems.
Earlier such
interpolation analogues for other
classical symmetric functions were obtained, see
\cite{Ok1}, \cite{OO1}, \cite{OO2},
\cite{Mo}, \cite{ORV}, \cite{KS}, \cite{OO3},
\cite{Kn}, \cite{Ok2}, \cite{Sa}.
The general scheme of such interpolation
for the symmetric polynomials was
considered in \cite{Ok3}.

\section{Pieri-type formula}\label{s-Pier}

\begin{defin}\label{d-stre}
For $\mu,\la\in DP$ we
will write $\mu\nearrow\la$,
if $|\mu|+1=|\la|$
and $\mu\subset\la$
(Definition~\ref{d-bord}).
\end{defin}

Recall that
$$
P_{(1); a} (x)= P_{(1)}^*(x)=P_{(1)}(x)=
\sum_j x_j.
$$

\begin{theor}\label{t-Pier}
As usual we assume $a_1=0$. Then we have
for an arbitrary $\mu\in DP$
\begin{equation}
P_{\mu; a}(P_{(1)}-\sum_{j=1}^{\ell(\mu)}
a_{\mu_j+1})=\sum_{\mu\nearrow\la} P_{\la; a}.
\label{f-Pier}
\end{equation}
\end{theor}

\begin{proof}
First we consider the case when all
$a_j$ are pairwise distinct.
Then (\ref{f-Pier}) may be easily
deduced from Theorem~\ref{t-char-I}.
In the general case we use the continuity
of the both sides of (\ref{f-Pier}).
\end{proof}

\begin{corol}\label{c-Pier}
For arbitrary $\mu\in DP$ we get
$$
P_\mu^*(P_{(1)}-|\mu|)=
\sum_{\mu\nearrow\la} P_\la^*.
$$
\end{corol}

\section{Dimensions of skew shifted
Young diagrams}
\label{s-dime}

\begin{defin}
Suppose $\mu$ and $\la$
are strict partitions and
$\mu\subset \la$ (Definition~\ref{d-bord}).
\emph{A standard shifted Young tableau
$T$ of shape $\la/\mu$} is
a bijection $T: D(\la)'\setminus
D(\mu)'\to\{1,2,\dots,|\la|-|\mu|\}$
such that $T(i,j)<T(i+1,j);\,\,
T(i,j)<T(i,j+1)$. The number of the standard
shifted tableaux of shape $\la/\mu$
is called \emph{a dimension of skew
shifted Young diagram $\la/\mu$}
and is denoted by $g_{\la/\mu}$.
\end{defin}

Now we give another description
of this dimension. First we define
a Schur graph (first considered  in \cite{Bo},
\cite{BO}).

\begin{defin}
The set of vertices of the Schur graph
is labelled by the set $DP$.
The Schur graph is an acyclic directed graph;
there is an edge directed from
$\mu\in DP$ to $\la\in DP$ if and only
if $\mu\nearrow\la$ (Definition~\ref{d-stre}).
\end{defin}

\begin{propo}
There exists a directed path
in the Schur
graph from $\mu\in DP$ to $\la\in DP$
if and only if $\mu\subset \la$.
In this case the number of these
directed paths is equal to
$g_{\la/\mu}$.
\end{propo}

\begin{proof}
The proof is trivial.
\end{proof}

An explicit formula
for $g_{\la/\varnothing}$
is well-known
(\cite[Ch.III, \S7, Examples]{Ma1}):
$$
g_{\la/\varnothing}=
\frac{|\la|!}{\prod_{k=1}^{\ell(\la)} \la_k!}
\prod_{i<j}\frac{\la_i-\la_j}{\la_i+\la_j}.
$$
In the next theorem we obtain an explicit
expression for $g_{\la/\mu}$.

\begin{theor}\label{t-dime}
Suppose $\mu,\la\in DP,\,\mu\subset\la$;
then
\begin{equation}
\label{f-dime}
g_{\la/\mu}=g_{\la/\varnothing}
\frac{P_\mu^*(\la_1,\dots,\la_{\ell(\la)})}
{(|\la|\fd |\mu|)}.
\end{equation}
\end{theor}

\begin{proof}
This theorem is proved in \cite{Iv1} by
a direct computation.
Also following \cite[\S9]{OO1} we may
deduce this theorem from the Pieri-type
formula (Corollary~\ref{c-Pier}).
\end{proof}

Earlier in \cite{OO1} A.~Okounkov
and G.~Olshanski studied the case of
the Young graph which is connected with the
Schur S-functions and the linear representations
of the symmetric groups.
Theorem~\ref{t-dime} may be viewed as
the projective analogue of \cite[Theorem 8.1]{OO1},
see also \cite{ORV}. The analogous result
for the more general case of
the Jack graph was obtained in
\cite[Section~5]{OO3}.

Now we rewrite the expression~(\ref{f-dime})
in the Pfaffian's form. In the ordinary case
the analogous formula
has a determinantal form, see, for example,
\cite[Proof of Proposition~1.2]{ORV}

\begin{theor}\label{t-dime-pfaf}
Suppose $\mu,\la\in DP,\,\mu\subset\la$.
Let $X(\la_1,\dots,\la_k)$ denote the skew-symmetric
$k\times k$ matrix
$$
\left(\frac{\la_i-\la_j}{\la_i+\la_j}\cdot
\frac1{\la_i!\la_j!}\right)_{1\le i,j\le k}
$$
and let $Y(\la_1,\dots,\la_k;\mu)$ denote the
$k\times \ell(\mu)$ matrix
$$
\left(\frac1{(\la_i-\mu_j)!}
\right)_{i\le k,\,j\le\ell(\mu)},
$$
where we suppose $\frac1{m!}=0$
if $m<0$. Put $s=\ell(\la)+\ell(\mu)$.
If $s$ is even then let $A_{\la/\mu}$
be the skew-symmetric $s\times s$ matrix
$$
\left(
\begin{array}{cc}
X(\la_1,\dots,\la_{\ell(\la)})&
Y(\la_1,\dots,\la_{\ell(\la)};\mu)\\
-Y(\la_1,\dots,\la_{\ell(\la)};\mu)^t&0
\end{array}
\right).
$$
If $s$ is odd then let $A_{\la/\mu}$
be the skew-symmetric $(s+1)\times (s+1)$
matrix
$$
\left(
\begin{array}{cc}
X(\la_1,\dots,\la_{\ell(\la)},0)&
Y(\la_1,\dots,\la_{\ell(\la)},0;\mu)\\
-Y(\la_1,\dots,\la_{\ell(\la)},0;\mu)^t&0
\end{array}
\right).
$$

Then
\begin{equation}
\label{f-Kiri}
g_{\la/\mu}=(|\la|-|\mu|)!\pf(A_{\la/\mu}).
\end{equation}
\end{theor}

\begin{proof}
The formula~(\ref{f-Kiri}) can be easily
deduced from Theorem~\ref{t-Nimm} and
Theorem~\ref{t-dime}.
\end{proof}

\section{Generating functions}
\label{s-gene}

It is known (see e.g. \cite[III, \S2, (2.10)]{Ma1})
that
$$
\sum_{r=0}^\infty Q_{(r)}(x_1, x_2,\dots)
y^r=\prod_{j=1}^\infty
\frac{1+x_j y}{1-x_j y}.
$$

\begin{defin}
Let $a=(a_n)_{n\ge 1}$
be a sequence of complex numbers,
then denote by $(\tau a)$
the result of shifting the sequence
$a$ to the left:
$$
(\tau a)_n= a_{n+1}.
$$
\end{defin}

Thus,
$$
(x|\tau a)^k=
\prod_{j=2}^{k+1}
(x-a_j).
$$

\begin{theor}
$$
\sum_{r=0}^\infty\frac{Q_{(r);a}(x_1,x_2,\dots)}
{(u|\tau a)^r}=
\prod_{j=1}^\infty\frac{u+x_j}{u-x_j}.
$$
\end{theor}

\begin{proof}
By our assumption (Definition~\ref{d-mult})
$a_1=0$. Then the following identity holds:
$$
1+2\sum_{r=1}^\infty
\frac{(x|a)^r}{(u|\tau a)^r}=
\frac{u+x}{u-x}.
$$
Then we may reason as in the proof of
\cite[III,\S2,(2.10)]{Ma1}.
\end{proof}

\begin{corol}
$$
\sum_{r=0}^{\infty}
\frac{Q^*_{(r)}(x_1,x_2,\dots)}{(u\fd r)}=
\prod_{i=1}^{\infty}\frac{u+1+x_i}{u+1-x_i}.
$$
\end{corol}

Next formula is due to A.~Borodin.

\begin{theor}\label{t-dvus}
If $0\le k\le l$ then put
$P_{(k,l); a}\equiv -P_{(l,k); a}$.
Then we have
\begin{multline}
\label{f-dvus}
4uv(u+v)\sum_{k,l\le 0}
\frac{P_{(k,l);a}(x_1,x_2,\dots)}
{(u|a)^{k+1}(v|a)^{l+1}}=\\
u\left(\prod_{j=1}^\infty
\frac{u+x_j}{u-x_j}+1\right)
\left(-\prod_{j=1}^\infty
\frac{v+x_j}{v-x_j}+1\right)-\\
v\left(-\prod_{j=1}^\infty
\frac{u+x_j}{u-x_j}+1\right)
\left(\prod_{j=1}^\infty
\frac{v+x_j}{v-x_j}+1\right).
\end{multline}
\end{theor}

\begin{proof}
The equality (\ref{f-dvus}) is equivalent to
following two equalities:
\begin{multline}
\label{f-dlin}
P_{(k+1,l);a}+P_{(k,l+1);a}+
(a_{k+1}+a_{l+1})P_{(k,l);a}=\\
P_{(k);a}P_{(l+1);a}-
P_{(k+1);a}P_{(l);a}+
(a_{l+1}-a_{k+1})P_{(k);a}P_{(l);a};
\end{multline}
\begin{equation}
\label{f-koro}
P_{(k+1);a}+P_{(k,1);a}+
a_{k+1}P_{(k);a}=P_{(k);a}P_{(1);a}.
\end{equation}
The equality
(\ref{f-koro}) is a particular case
of Theorem~\ref{t-Pier}.
Let us prove (\ref{f-dlin}).
From \cite[Ch.III, \S5, (5.7)]{Ma1}
we get the following identity
for the ordinary Schur P-functions
$$
P_{(k+1)}P_{(l)}=
P_{(k+1,l)}+P_{(k+1+l)}+
2\sum_{j=1}^{l-1}P_{(k+1+j,l-j)}.
$$
Consequently we have
\begin{equation}
\label{f-simp}
P_{(k+1,l)}+P_{(k,l+1)}=
P_{(k)}P_{(l+1)}-
P_{(k+1)}P_{(l)}.
\end{equation}

First we suppose that the numbers
$a_j$ are pairwise distinct;
then from (\ref{f-simp})
and Theorem~\ref{t-char-I} we get
that
\begin{multline}
\label{f-cons}
P_{(k+1,l);a}+P_{(k,l+1);a}+
\alpha P_{(k,l);a}=\\
P_{(k);a}P_{(l+1);a}-
P_{(k+1);a}P_{(l);a}+
(a_{l+1}-a_{k+1})P_{(k);a}P_{(l);a}
\end{multline}
for some $\alpha$.
Evaluating both sides of (\ref{f-cons})
at the point $(a_{k+1}, a_{l+1})$
we obtain that $\alpha=a_{k+1}+a_{l+1}$.
Thus (\ref{f-dlin}) is proved for sequences
$(a_j)$ without repetitions.
Using the continuity argument we have that
(\ref{f-dlin}) holds for all sequences
$(a_j)$. This concludes the proof.
\end{proof}

One may compare Theorem~\ref{t-dvus}
with its analog
for the multiparameter
supersymmetric Schur functions
(\cite[Proposition~7.1]{ORV}).

\section{Giambelli-Schur-type
formula}\label{s-Giam}

\begin{theor}
\label{t-Giam}
As in Theorem~\ref{t-dvus},
we put $Q_{(k,l); a}\equiv
-Q_{(l,k); a}$ for arbitrary
integers $l\ge k\ge 0$.
Also we will use the notation
$\la_{\ell(\la)+1}$, certainly,
$\la_{\ell(\la)+1}=0$. If $\la$ is
a strict partition then we have
$$
Q_{\la;a}=\pf
\left( (Q_{(\la_i,\la_j);
a})_{1\le i,j\le 2[\frac{\ell(\la)+1}2]}
\right),
$$
where $[x]$ stands for the integral
part of $x$.
\end{theor}

\begin{proof}
We use Theorem~\ref{t-comb} and the
method of Stembridge \cite{St2}
with simple modifications.
\end{proof}

I.~Schur (\cite{Sc}) used the particular case
of this formula $(a_j\equiv 0)$
as the definition of the Q-functions.
This formula may be viewed as the projective
analogue of the Giambelli formula for
the ordinary s-functions
\cite[Ch.I, \S3, Examples]{Ma1}.
We may compare Theorem~\ref{t-Giam}
with its determinantal analogue
for the multiparameter
supersymmetric Schur functions
\cite{ORV}.

\section{Transition coefficients}
\label{s-tran}

In this paragraph we use the method of
\cite[\S2, \S7]{ORV}.
From \cite[Lemma 2.5]{ORV} we get the
following identity for arbitrary
sequences $(a_n)_{n\ge 1}$ and $(b_n)_{n\ge 1}$
\begin{multline}
\label{f-hsup}
\frac1{(u-b_1)\dots (u-b_{r'+1})}=\\
\sum_{r=r'}^\infty\frac
{h_{r-r'}(b_1,\dots,b_{r'+1};-a_1,\dots,-a_r)}
{(u-a_1)\dots(u-a_{r+1})},
\end{multline}
where $h_0=1,h_1,h_2,\dots$ denote the
conventional complete homogenous functions
in the super realization of the algebra
$\Lambda$ of the symmetric functions
(\cite[Ch.I, \S5, Examples]{Ma1}).
From now in this
paragraph suppose that the sequences
$(a_n)_{n\ge 1}$ and $(b_n){n\ge 1}$
are fixed and $a_1=b_1=0$. Then
\begin{multline}
h_k(b_1,b_2,\dots,b_{r'+1};
-a_1,-a_2,\dots,-a_r)=\\
h_k(b_2,\dots,b_{r'+1};
-a_2,\dots,-a_r).
\end{multline}
Denote by $d_{rr'}$ the value
$h_k(b_2,\dots,b_{r'+1};
-a_2,\dots,-a_r)$, where
$r\ge r'$. If $r<r'$ then put
$d_{rr'}=0$.

\begin{propo}
\label{p-tran}
Let $r,s\ge 1$. Then
$$
P_{(r,s);a}=\sum_{r'=1}^r
\sum_{s'=1}^s
d_{rr'}(a,b) d_{ss'}(a,b)
P_{(r',s');b}.
$$
\end{propo}

\begin{proof}
The proof is based on
(\ref{f-hsup}) and Theorem~\ref{t-dvus}.
\end{proof}

\begin{theor}\label{t-tran}
Suppose $\mu$ is a strict partition;
then
$$
Q_{\mu;a}=\sum_{\nu\in DP,
\ell(\nu)=\ell(\mu), \nu\subset\mu}
d_{\mu\nu} Q_{\nu;b},
$$
where
$$
d_{\mu\nu}=\det
\left((d_{\mu_i\nu_j})_{1\le i,j\le\ell(\mu)}
\right).
$$
\end{theor}

\begin{proof}
We use Theorem~\ref{t-Giam} and
Proposition~\ref{p-tran}.
\end{proof}

Theorem~\ref{t-tran} is the projective
analogue of \cite[Theorem~7.3]{ORV}.

\end{document}